# Asymptotics of the number of primes and sums of functions of primes in a subset of natural series

Victor Volfson

ABSTRACT The paper solves the problems of determining the asymptotics of the number of primes and the sums of functions of primes in a subset of the natural series that satisfies the conditions that the asymptotic density of the number of primes in this subset is constant and not equal to zero.

Assertions about asymptotics of the sums of functions of primes in a subset of the natural series satisfying the indicated conditions are proved. Necessary and sufficient conditions for the existence of these asymptotics are also proved.

## 1. INTRODUCTION

Dirichlet proved in 1837 that any arithmetic progression whose initial term is relatively prime to the difference contains an infinite number of primes.

There are only conjectures for an integer polynomial of degree greater than one when it contains an infinite number of primes. One of them is Bunyakovskii's conjecture [1]. If $f(n)$ is an integer irreducible polynomial and $\mathrm{a}$ is the greatest common divisor of all its values, then the integer polynomial $f(n)/\mathrm{a}$ takes an infinite number of prime values.

Conjecture that the polynomial $f(n)=n^2+1$ takes on an infinite number of prime values is a particular conjecture of Bunyakovsky (Landau's 4th problem).

Asymptotics of the number of primes in the natural series and the arithmetic progression under the Dirichlet conditions was first proved in [2].

Asymptotics of the sums of functions of prime numbers in natural series and arithmetic progression, under the Dirichlet conditions, was proved in [3].

There is only the Bateman-Horn conjecture [4] about asymptotics of the number of primes in an integer polynomial of degree higher than the first under the conditions of the Bunyakovskii conjecture. Asymptotics of the number of natural numbers on the interval $[2,\mathrm{x})$ for which the value of the degree $\mathrm{h}$ polynomial $f$ is a prime number is determined by the formula:

$$\pi(\mathrm{f},\mathrm{x}) \sim \frac{c(f)}{h}\int_2^x \frac{dt}{\ln t}, \qquad (1.1)$$

where $c(f)=\prod_p \frac{1-\omega(p)/\mathrm{p}}{1-1/p}$, $\omega(p)$ is the number of comparison solutions $f\equiv 0(\mathrm{mod}\, p)$.

There are not even conjectures about the asymptotics of the sums of functions of prime numbers that are in an integer polynomial of degree higher than the first.

Therefore, the problem of determining the asymptotics of the number of primes and the sums of functions of primes in a subset of the natural series that satisfies certain requirements is topical.

Let's start with definitions. An integer sequence is a strictly increasing sequence of natural numbers. An example of an integer sequence is the sequence of prime numbers: 2, 3, 5, 7,…

Denote any integer sequence by the letter $A$, the number of members of this integer sequence at $n \le x$ - $\pi(\mathrm{A},\mathrm{x})$, i.e.

$$\pi(\mathrm{A},\mathrm{x}) = |\{n \le x : n \in A\}| . \tag{1.2}$$

Having in mind (1.2), the density of an integer sequence $A$ at value $n \le x$ is determined by the formula:

$$\mathrm{d}(\mathrm{A},\mathrm{x}) = \pi(\mathrm{A},\mathrm{x}) / \mathrm{x} = |\{n \le x : n \in A\}| / x . \tag{1.3}$$

Based on the monotonicity of an integer sequence, the density of an integer sequence (1.3) is a probability measure or simply a probability.

Let's denote the asymptotic or natural density of an integer sequence $A$ as:

$$d^*(A) = \lim_{x\to\infty} \pi(\mathrm{A},\mathrm{x}) / \mathrm{x} = \lim_{x\to\infty} |\{\mathrm{n} \le x : n \in A\}| / x . \tag{1.4}$$

The asymptotic density of any integer sequence with a finite number of terms is equal to zero, so the asymptotic density of an integer sequence is not a measure on the natural series. An integer sequence with an asymptotic density does not form a sigma-algebra. Moreover, the asymptotic density is not countably additive [5].

Considering the above, the asymptotic density of an integer sequence is not a probability measure, but due to (1.3), (1.4) it is its limit.

Since an integer sequence is a subset of the natural series, then, based on (1.4), the following holds:

$$d^*(A) \le 1 . \tag{1.5}$$

We will consider a sequence of prime numbers - $P$ with members $\mathrm{p} \in \mathrm{P}$ as an integer sequence in this work.

We will omit the integer sequence notation $A$ when we talk about sequences of primes.

Based on the asymptotic law of primes [2], the density of primes, taking into account (1.3), is determined by the formula:

$$d(\mathrm{x}) = \pi(\mathrm{x}) / \mathrm{x} = |\{\mathrm{n} \le \mathrm{x} : n \in P\}| / x = \frac{1}{\ln x}(1 + o(1)) \,. \qquad (1.6)$$

This density is the probability of a randomly chosen natural number not exceeding $x$, to be a prime.

## 2.ASYMPTOTICS OF THE NUMBER OF PRIMES IN A SUBSET OF THE NATURAL SERIES

Now we will consider a subset of the natural series $f$, which is also an integer sequence. Let us denote the number of natural numbers, belonging to the intersection $f \cap A$, that do not exceed the value $x$ - $\pi(\mathrm{A,f,x})$:

$$\pi(\mathrm{A}, f, \mathrm{x}) = |\{\mathrm{n} \le \mathrm{x} : \mathrm{n} \in f \bigcap A\}| \,. \qquad (2.1)$$

Let us denote the density of natural numbers that do not exceed the value $x$, belonging to the intersection $f \cap A$:

$$d(\mathrm{A}, f, x) = \pi(A, f, x) / \pi(A, x) = |\{n \le \mathrm{x} : \mathrm{n} \in f \cap A\}| / |\{n \le \mathrm{x} : n \in \mathrm{A}\} \,. \qquad (2.2)$$

Let us denote the asymptotic density of natural numbers, belonging to the intersection $f \cap A$:

$$\mathrm{d}^*(\mathrm{A}, f) = \lim_{x\to\infty} \pi(A, f, x) / \pi(A, x) = \lim_{x\to\infty} |\{n \le \mathrm{x} : n \in f \cap A\}| / |\{n \le \mathrm{x} : n \in A\}| . \qquad (2.3)$$

Then, based on (2.1) and (2.2), the number of natural numbers that do not exceed the value $x$, belonging to the intersection $f \cap A$, is equal to:

$$\pi(A, f, x) = d(\mathrm{A}, f, x)\pi(A, x) \,. \qquad (2.4)$$

If $\mathrm{d}^*(\mathrm{A}, f)$ is a constant not equal to zero, then based on (2.3) we get:

$$\lim_{x\to\infty} \frac{\pi(A, f, x)}{\mathrm{d}^*(\mathrm{A}, f)\pi(A, x)} = 1 \,. \qquad (2.5)$$

Having in mind (1.5) we have $0 < \mathrm{d}^*(\mathrm{A}, f) \le 1$.

Based on (2.5), the asymptotic of the number of natural numbers that do not exceed the value $x$ , belonging to the intersection $f \cap A$, is equal to:

$$\pi(A, f, x) \sim \mathrm{d}^*(\mathrm{A}, f)\pi(A, x) \text{ или } \pi(A, f, x) = \mathrm{d}^*(\mathrm{A}, f)\pi(A, x)(1+o(1)) \ (\mathrm{x} \to \infty) . (2.6)$$

In a particular case, for the sequence of primes, based on (2.6), the asymptotic of the number of natural numbers that do not exceed the value $x$, belonging to the intersection $P \cap f$, is equal to:

$$\pi(f, x) = \mathrm{d}^*(f)\pi(x)(1+o(1)) = \mathrm{d}^*(f)\frac{x(1+o(1))}{\ln x} \ (\mathrm{x} \to \infty) . \qquad (2.7)$$

Therefore, determining the asymptotic of the number of natural numbers that do not exceed the value $x$, belonging to the intersection $P \cap f$, in the case when $\mathrm{d}^*(f)$ is a constant and not equal to zero, is reduced to finding the asymptotic density $\mathrm{d}^*(f)$ using formula (2.3).

Now let's consider examples of finding the asymptotic of natural numbers that do not exceed the value $x$ , belonging to the intersection $P \cap f$, in the case when $\mathrm{d}^*(f)$ is a constant not equal to zero.

First, let's look at trivial examples to clarify the essence, and then consider a not-quite-trivial example.

As a first example, we will consider arithmetic progressions with a denominator - 3: $3n, 3n+1, 3n+2$ .

The subset $f = 3n$ does not take prime values, so the intersection $P \cap f$, in this case is empty and $\mathrm{d}^*(f) = 0$, so the specified subset does not satisfy the above requirements.

Arithmetic progressions: $3n+1, 3n+2$ take prime values in equal parts, so the asymptotic densities $\mathrm{d}^*(f)$, in these subsets, are:

$$d^*(3n+1) = d^*(3n+2) = 0{,}5 . \qquad (2.8)$$

and satisfies the above requirements.

Based on (2.7) and (2.8), asymptotic of natural numbers, belonging to the intersection $P \cap f$, for subsets $3n+1, 3n+2$, respectively, are:

$$\pi(3n+1,\mathrm{x}) = \pi(3n+2,\mathrm{x}) = \frac{0,5x(1+o(1))}{\ln x} \ (\mathrm{x} \to \infty). \qquad (2.9)$$

In general, there are $\varphi(\mathrm{k})$ arithmetic progressions $kn+l(k,l=1)$ that take prime numbers in equal parts, so the asymptotic densities $\mathrm{d}^*(f)$ in given subsets of the natural series are equal $1/\varphi(\mathrm{k})$ and satisfy the above requirements.

Having in mind (2.7), asymptotic of natural numbers, belonging to the intersection $P \cap f$, for these subsets, are equal to:

$$\pi(kn+l,x) = \frac{x(1+o(1))}{\varphi(\mathrm{k})\ln x} \ (\mathrm{x} \to \infty). \qquad (2.10)$$

Now let's define the asymptotic of natural numbers belonging to the intersection $P \cap f$, for the subset $f$ - numbers free from squares.

Natural numbers that are in the subset of square-free numbers take prime values, and all of them, so the asymptotic density of natural numbers belonging to the intersection $P \cap f$, for this subset is equal to 1 and this subset satisfies the above requirements.

Based on (2.7), the asymptotic of natural numbers, belonging to the intersection $P \cap f$, for the subset of square-free numbers is:

$$\frac{x(1+o(1))}{\ln x} \ (\mathrm{x} \to \infty). \qquad (2.11)$$

Now there will be a not quite trivial example. The natural numbers that divide an arbitrary prime number - $p$ are at an equal distance - $p$, so the asymptotic density of such numbers in the natural series is:

$$d^*(p) = 1/p. \qquad (2.12)$$

By analogy with probability, one can introduce the concept of independence of the asymptotic densities of given numbers.

Let there are prime numbers $p_i, p_j (i \neq j)$.

Then, using (2.12), one can show the independence of the indicated asymptotic densities:

$$\mathrm{d}(p_i p_j) = \frac{1}{p_i p_j} = \frac{1}{p_i}\frac{1}{p_j} = d(p_i)d(p_j)\,. \qquad (2.13)$$

However, the asymptotic densities are not probabilities, so the product $\prod_{p\le\sqrt{x}}(1-1/p)$ is not the probability of a randomly chosen natural number not exceeding $x$, to be a prime.

On the other hand, the Mertens theorem [5] is true:

$$\frac{e^{-\gamma}}{\ln y} = \prod_{p\le y}(1-1/p)(1+o(1))\,. \qquad (2.14)$$

Let us make a change of variables in (2.14) and write it in the form:

$$\frac{1}{\ln x}(1+o(1)) = 0{,}5e^{\gamma}\prod_{p\le\sqrt{x}}(1-1/p)\,. \qquad (2.15)$$

Based on (1.6), there is a probability (on the left in (2.15)) of a randomly chosen natural number not exceeding $x$, to be a prime.

Now let's consider the asymptotic of the number of natural numbers for which $f = n^2+1 \in P$.

Based on (1.1) we have in this case:

$$\pi(f,\mathrm{x}) \sim \frac{c(f)}{2}\int_2^x \frac{dt}{\ln t}, \qquad (2.16)$$

where $$\frac{c(f)}{2} = \frac{\prod_{p\ge 2}(1-\omega(p)/p)}{\prod_{p>2}(1-1/\mathrm{p})} = 0{,}6864067...\,, \qquad (2.17)$$

$\omega(p)$ is the number of comparison solutions $x^2+1 \equiv 0(\mathrm{mod}\,p)$.

We write (2.17) in the form:

$$\frac{c(f)}{2} = \lim_{x\to\infty}\frac{0{,}5e^{\gamma}\prod_{p\le\sqrt{x}}(1-\omega(p)/p)}{0{,}5e^{\gamma}\prod_{p\le\sqrt{x}}(1-1/\mathrm{p})}\,. \qquad (2.18)$$

Based on (2.15), the probability of a randomly chosen natural number not exceeding $x$, to be a prime is in the denominator of (2.18). It is equal to the number of prime numbers not exceeding $x$, divided by $x$. The numerator (2.18) contains the probability of a randomly chosen natural number not exceeding $x$, for which the value $n^2+1$ is a prime. It is equal to the number of prime numbers not exceeding $x$, for which the value $n^2+1$ is a prime, divided by $x$.

Hence:

$$\frac{c(f)}{2} = \lim_{x\to\infty} \frac{\frac{|\{\mathrm{n} \le x : n^2+1 \in P\}|}{x}}{\frac{|\{n \le x : n \in P\}|}{x}} = \lim_{x\to\infty} \frac{|\{\mathrm{n} \le x : n^2+1 \in P\}|}{|\{n \le x : n \in P\}|} . \tag{2.19}$$

Therefore, based on (2.19), $\frac{c(f)}{2}$ is the asymptotic density $d^*(f)$ for the subset $n^2+1$, i.e. $d^*(f) = \frac{c(f)}{2}$.

Since $\frac{c(f)}{2}$ is a constant not equal to 0, i.e. satisfies our conditions, then the asymptotic of natural numbers belonging to the intersection $P \bigcap f$, for the subset $n^2+1$, based on (2.7), is indeed determined by the formula (2.16).

Based on (2.6) and the law of prime numbers [2], the asymptotic of natural numbers that do not exceed $x$ and belonging to the intersection $P \bigcap f$, in the case, when $\mathrm{d}^*(f)$ is a constant not equal to zero, is determined by the formula:

$$\pi(f, \mathrm{x}) = \mathrm{d}^*(f)\int_2^x \frac{dt}{t} + O(\frac{x}{e^{c\sqrt{\ln x}}}) , \tag{2.20}$$

where $\mathrm{c} > 0$ is a constant.

Using (2.6) and integrating (2.20) by parts, we obtain another formula for determining the given asymptotic:

$$\pi(f, \mathrm{x}) = \mathrm{d}^*(f)\frac{x}{\ln x} + O(\frac{x}{\ln^2 x}) . \tag{2.21}$$

If the Riemann conjecture is true, then, based on (2.6), this asymptotic behavior is determined by the formula:

$$\pi(f,\mathrm{x}) = \mathrm{d}^*(f)\int_2^x \frac{dt}{t} + O(\mathrm{x}^{1/2}\ln\mathrm{x}). \qquad (2.22)$$

In Chapter 3, based on these asymptotics, a general formula will be obtained for the asymptotic estimate of the sums of functions of prime numbers that are on a subset of the natural series (that satisfies the specified requirements), and then estimates of asymptotics for various functions will be obtained.

In Chapter 4, necessary and sufficient conditions for the existence of these asymptotics will be proved.

3. ASYMPTOTICS OF SUMS OF FUNCTIONS OF PRIME NUMBERS IN A SUBSET OF THE NATURAL SERIES

Let $a_m$ is a sequence of real numbers and $g(m)$ is a continuously differentiable function on the ray $[2,x)$, then the Abel summation formula is true [5]:

$$\sum_{p\le x, p\in f} g(p) = \sum_{m=2}^{x} a_m g(m) = A(\mathrm{x})g(x) - \int_2^x A(t)g'(t)dt, \qquad (3.1)$$

where $A(x) = \sum_{m=1}^{x} a_m$ .

Let's look at the distribution of primes in a subset $f$ of the natural series. Let $a_m = 1$ if the value $f(m)$ is a prime number, otherwise $a_m = 0$. Then $A(x) = \pi(f,x)$, where $\pi(f,x)$ is the number of prime numbers in the subset $f$ of natural numbers, that do not exceed $x$.

Substitute $A(x) = \pi(f,x)$ in formula (3.1) and we get:

$$\sum_{p\le x, p\in f} g(p) = \sum_{m=2}^{x} a_m g(m) = \pi(f,x)g(x) - \int_2^x \pi(f,\mathrm{t})g'(t)dt. \qquad (3.2)$$

We will use (3.2) to prove a number of assertions.

We use the asymptotic estimate (2.21) to determine the asymptotic behavior of the sums of prime functions that are in a subset $f$ of the natural series.

Assertion 3.1

If $g(t)$ is a monotonic function and has a continuous derivative on the ray $[2,x)$, and if $d^*(f)$ is a constant not equal to 0, then:

$$\sum_{p\le x, p\in f} g(p) = \frac{d^*(f)xg(x)}{\ln(x)} - d^*(f)\int_2^x \frac{tg'(t)dt}{\ln(t)} + O(\frac{x\,|\,g(x)\,|}{\ln^2(x)}) + O(\int_2^x \frac{t\,|\,g'(t)\,|\,dt}{\ln^2(t)}) . \quad (3.3)$$

To prove this, we substitute (2.21) into (3.2) and obtain (3.3).

Let us give examples of using formula (3.3):

1. $\sum_{p\le x, p\in f} 1 = \frac{d^*(f)x}{\ln(x)} + O(\frac{x}{\log^2(x)})$, which corresponds to formula (2.21).

2. $\sum_{p\le x, p\in a} \ln(p) = d^*(\mathrm{f})\,\mathrm{x} + O(x/\ln(x))$. (3.4)

3. $\sum_{p\le x, p\in f} 1/p = d^*(\mathrm{f})\ln\ln(x) + O(1)$. (3.5)

Now we use a more precise asymptotic estimate (2.22).

Assertion 3.2

If the function $f$ is monotonic and has a continuous derivative on the ray $[2,x)$, then if $d^*(f)$ is a constant not equal to 0, then:

$$\sum_{p\le x, p\in f} g(p) = d^*(f)\int_2^x \frac{g(t)dt}{\ln(t)} + O(\frac{|\,g(x)\,|\,x}{e^{c\ln^{1/2}(x)}}) + O(\int_2^x \frac{t\,|\,g'(t)\,|\,dt}{e^{c\ln^{1/2}(t)}}) ,$$

where $c$ is a constant.

Proof

We substitute (2.20) into (3.2) and get:

$$\sum_{p\le x, p\in f} g(p) = d^*(f)g(x)\int_2^x \frac{dt}{\ln(t)} + O(\frac{|g(x)|x}{e^{c\ln^{1/2}(x)}}) - d^*(f)\int_2^x (\int_2^t \frac{du}{\ln(u)})g'(t)dt + O(\int_2^x \frac{t|g'(t)|dt}{e^{c\ln^{1/2}(t)}}) .(3.6)$$

We use the integration by parts method:

$$d^*(f)\int_2^x (\int_2^t \frac{du}{\ln(u)})g'(t)dt = d^*(f)g(x)\int_2^x \frac{du}{\ln(u)} - d^*(f)\int_2^x \frac{g(t)dt}{\ln(t)} . \qquad (3.7)$$

Substitute (3.7) into (3.2) and get:

$$\sum_{p\le x, p\in f} g(p) = d^*(f)\int_2^x \frac{g(t)dt}{\ln(t)} + O(\frac{|g(x)|x}{e^{c\ln^{1/2}(x)}}) + O(\int_2^x \frac{t|g'(t)|dt}{e^{c\ln^{1/2}(t)}}) , \qquad (3.8)$$

which corresponds to assertion 3.2.

Let's look at examples of use (3.8):

1. $\sum_{p\le x, p\in f} 1 = d^*(f)\int_2^x \frac{dt}{\ln(t)} + O(\frac{x}{e^{c\ln^{1/2}(x)}})$, which corresponds to (2.20).

2. $\sum_{p\le x, p\in f} \ln(p) = d^*(f)\int_2^x \frac{\ln(t)dt}{\ln(t)} + O(\frac{\ln(x)x}{e^{c\ln^{1/2}(x)}}) + O(\int_2^x \frac{dt}{e^{c\ln^{1/2}(t)}}) = d^*(f)x + O(\frac{x\ln(x)}{e^{c\ln^{1/2}(x)}})$ . (3.9)

Compare (3.9) and (3.4).

3. $\sum_{p\le x, p\in f} \frac{\ln(p)}{p} = d^*(f)\int_2^x \frac{\ln(t)dt}{t\ln(t)} + O(\frac{\ln(x)x}{xe^{c\ln^{1/2}(x)}}) + O(\int_2^x \frac{tdt}{t^2e^{c\ln^{1/2}(t)}}) + O(\int_2^x \frac{t\ln(t)dt}{t^2e^{c\ln^{1/2}(t)}}) = d^*(f)\ln x + O(1)$

4. $\sum_{p\le x, p\in f} p^\alpha = d^*(f)\int_2^x \frac{t^\alpha dt}{\ln(t)} + O(\frac{x^{\alpha+1}}{e^{c\ln^{1/2}(x)}}) + O(\int_2^x \frac{t^\alpha dt}{e^{c\ln^{1/2}(t)}}) = d^*(f)\int_2^x \frac{t^\alpha dt}{\ln(t)} + O(\frac{x^{\alpha+1}}{e^{c\ln^{1/2}(x)}})$ ,

where $\alpha > -1$.

Assertion 3.3

If the function $f$ is monotonic and has a continuous derivative on the ray $[2,x)$, then if the generalized Riemann conjecture is true and if $d^*(f)$ is a constant not equal to 0, then:

$$\sum_{p\le x, p\in f} g(p) = d^*(f)\int_2^x \frac{g(t)dt}{\ln(t)} + O(|g(x)|x^{1/2}\ln(x)) + O(\int_2^x |g'(t)|t^{1/2}\ln(t)dt) . \qquad (3.10)$$

The proof is carried out similarly to assertion 3.2, only the value (2.22) is substituted into (3.2).

The asymptotic estimate for the remainder term (3.10) is the best. This can be verified by comparing the remainder terms in Assertion 3.3 with Assertions 3.1 and 3.2.

Let us show that the main terms of the asymptotics (3.3) and (3.8) coincide, and the differences are in the remainder terms.

We transform the term in formula (3.8):

$$\int_2^x \frac{g(t)dt}{\ln(t)} = \frac{g(x)x}{\ln(x)} - \frac{2f(2)}{\ln 2)} - \int_2^x \frac{tg'(t)dt}{\ln(t)} + \int_2^x \frac{g(t)dt}{\ln(t)} . \quad (3.11)$$

We substitute (3.11) into (3.8):

$$\sum_{p\le x, p\in f} g(p) = \frac{d^*(f)g(x)x}{\ln(x)} + O(1) - d^*(f)\int_2^n \frac{tg'(t)dt}{\ln(t)} + d^*(f)\int_2^n \frac{g(t)dt}{\ln(t)} + O(\frac{|g(x)x|}{e^{c\ln^{1/2}(x)}}) + O(\int_2^x \frac{t|g'(t)|dt}{e^{c\ln^{1/2}(t)}}) . \quad (3.12)$$

Let us compare expression (3.12) with the asymptotic expression (3.3). The first three terms are the same. The discrepancy is only in the remainder terms. The remainder terms of expression (3.12) give a more precise estimate.

## 4. NECESSARY AND SUFFICIENT CONDITIONS FOR THE EXISTENCE OF THESE ASYMPTOTICS

Let $a_m = 1$ if $f(\text{m})$ is a prime number, and $a_m = 0$ otherwise. Let $b_1 = 0, b_m = \frac{f(m)}{\ln m}$, where $m$ is a natural number.

Let's denote: $A(\text{n}) = \sum_{m\le n} a_m = \sum_{m\le n} a_m g(m)$ and $B(n) = \sum_{m\le n} b_m = \sum_{m\le n} b_m g(m)$.

It is required that:

$$\lim_{n\to\infty} \frac{A(n)}{B(n)} = 1 . \quad (4.1)$$

We write (4.1) in the form:

$$\lim_{n\to\infty}\frac{\sum_{k=1}^{n} a_k g(k)}{\sum_{k=1}^{n} b_k g(k)} = 1. \quad (4.2)$$

Assertion 4.1

When the conditions are met:

1. $\lim_{n\to\infty}\frac{\int_1^n B(t)g'(t)dt}{B(n)g(n)}$ is not equal to 1.

2. $g(x)$ - monotonic and $g'(x) \neq 0$.

3. $\lim_{n\to\infty}\int_1^n B(t)\,\mathrm{g}'(t)dt = \pm\infty$.

Then it is performed:

$$\lim_{n\to\infty}\frac{\sum_{k=1}^{n} a_k g(k)}{\sum_{k=1}^{n} b_k g(k)} = 1$$

Proof

Having in mind the Abel summation formula:

$$\sum_{m=1}^{n} a_m g(m) = A(n)g(n) - \int_1^n A(t)g'(t)dt\,,\ \sum_{m=1}^{n} b_m\, \mathrm{g}(m) = B(n)g(n) - \int_1^n B(t)g'(t)dt\,. \quad (4.3)$$

Substituting expressions (3.3) into (3.2) we get:

$$\frac{\sum_{k=1}^{n} a_k g(k)}{\sum_{k=1}^{n} b_k g(k)} = \frac{A(n)g(n) - \int_1^n A(t)g'(t)dt}{B(n)g(n) - \int_1^n B(t)g'(t)dt} = \frac{A(n)}{B(n)}\cdot\frac{1-\frac{\int_1^n A(t)g'(t)dt}{A(n)g(n)}}{1-\frac{\int_1^n B(t)g'(t)dt}{B(n)g(n)}}\,. \quad (4.4)$$

Let the above sufficient conditions is satisfied:

1. $\lim_{n\to\infty}\frac{\int_1^n B(t)g'(t)dt}{B(n)g(n)}$ is not equal to 1.

2. Let $g(x)$ - monotonous and $g'(x) \neq 0$.

3. $\lim_{n\to\infty} \int_1^n B(t)g'(t)dt = \pm\infty$.

Let us show that under the above conditions:

$$\lim_{n\to\infty} \frac{\int_1^n A(t)g'(t)dt}{A(n)g(n)} = \lim_{n\to\infty} \frac{\int_1^n B(t)g'(t)dt}{B(n)g(n)} \quad . \tag{4.5}$$

Since $\lim_{n\to\infty} \frac{A(n)}{B(n)} = 1$, then using L'Hopital's rule, we get:

$$\lim_{n\to\infty} \frac{\dfrac{\int_1^n A(t)g'(t)dt}{A(n)g(n)}}{\dfrac{\int_1^n B(t)g'(t)dt}{B(n)g(n)}} = \lim_{n\to\infty} \frac{B(n)}{A(n)} \cdot \lim_{n\to\infty} \frac{\int_1^n A(t)g'(t)dt}{\int_1^n B(t)g'(t)dt} = \lim_{n\to\infty} \frac{A(n)g'(n)}{B(n)g'(n)} = 1, \tag{4.6}$$

which corresponds to (4.5).

Then, based on (4.1) and (4.6), we obtain:

$$\lim_{n\to\infty} \frac{\sum_{k=1}^{n} a_k g(k)}{\sum_{k=1}^{n} b_k g(k)} = \lim_{n\to\infty} \frac{A(n)}{B(n)} \cdot \frac{1 - \dfrac{\int_1^n A(t)g'(t)dt}{A(n)g(n)}}{1 - \dfrac{\int_1^n B(t)g'(t)dt}{B(n)g(n)}} = 1,$$

which is consistent with the assertion.

Corollary 4.2

Conditions (1), (3) in Assertion 4.1 correspond to:

1. $\lim_{n\to\infty} \frac{\int_2^n \frac{tg'(t)}{\ln(t)}dt}{\frac{ng(n)}{\ln(n)}}$ is not equal to 1.

3. $\lim_{n\to\infty} \int_2^n \frac{tg'(t)}{\ln(t)}dt = \pm\infty$.

Proof

Let's find the value:

$$B(n) = \sum_{m \le n} b_m = d^*(f) \sum_{m \le n} \frac{1}{\ln m} = \frac{d^*(f) n}{\ln(\mathrm{n})}(1 + o(1)) . \qquad (4.7)$$

We substitute (4.7) into conditions 3 of assertion 4.1 and obtain:

$$\lim_{n\to\infty} \int_1^n B(t) g'(t) dt = d^*(f) \int_2^n \frac{t g'(t) dt}{\ln t}(1 + o(1)) . \qquad (4.8)$$

Based on (4.8), condition (3) is satisfied if $\lim_{n\to\infty} \int_2^n \frac{t g'(t) dt}{\ln t} = \pm\infty$ .

Now we substitute (4.7) into condition 1 of assertion 4.1 and get:

$$\lim_{n\to\infty} \frac{\int_1^n B(t) g'(t) dt}{B(n) g(n)} = \lim_{n\to\infty} \frac{d^*(f) \int_2^n \frac{t g'(t) dt}{\ln t}(1 + o(1))}{\frac{d^*(f) n g(n)}{\ln(\mathrm{n})}(1 + o(1))} = \lim_{n\to\infty} \frac{\int_2^n \frac{t g'(t) dt}{\ln t}}{\frac{n g(n)}{\ln(\mathrm{n})}}$$ is not equal to 1.

The conditions of Assertion 4.1 and Corollary 4.2 are sufficient for the above asymptotics to hold.

Now we consider a particular case of Assertion 4.1 and Corollary 4.2.

Assertion 4.3

Let $g(n)$ - increase, i.e. $\mathrm{g}'(n) > 0$, $\lim_{n\to\infty} g(n) = \infty$, $g'(n)$ is a continuous function and $\lim_{n\to\infty} \frac{g(n)}{\ln(n) g'(n)} \ne 0$. Then all the conditions of Assertion 4.1 and Corollary 4.2 are satisfied.

Proof

Let us check the 1st condition of Corollary 4.2 using L'Hopital:

$$\lim_{n\to\infty} \frac{\int_2^n \frac{t g'(t)}{\ln(t)} dt}{\frac{n g(n)}{\ln(n)}} = \lim_{n\to\infty} \frac{\frac{n g'(n)}{\ln(n)}}{(\frac{n g(n)}{\ln(n)})'} = \lim_{n\to\infty} \frac{\frac{n f'(n)}{\log(n)}}{\frac{n g'(n)}{\ln(n)} + (\frac{n}{\ln(n)})' g(n)} = \lim_{n\to\infty} \frac{1}{1 + \lim_{n\to\infty} (\frac{g(n)}{n g'(n)} - \frac{g(n)}{\ln(n) g'(n)})} \ne 1 . (4.9)$$

It is required $\lim_{n\to\infty}\frac{g(n)}{\ln(n)g'(n)}\neq 0$ that (4.9) was satisfied.

Let us check the second condition of Assertion 4.1.

$g(n)$ is a monotone function and $g'(n)\neq 0$ by the condition of Assertion 4.3.

Let us check the 3rd condition of Corollary 4.2 and show that $\lim_{n\to\infty}\int_2^n \frac{tg'(t)}{\ln(t)}dt=\pm\infty$.

$\int_2^n \frac{tg'(t)}{\ln(t)}dt \geq \int_2^n \frac{2g'(t)}{\ln(2)}dt = \frac{2}{\ln(2)}\int_2^n g'(t)dt = \frac{2}{\ln(2)}(g(n)-f(2))$ - increases indefinitely, as $g(n)$.

Let us show that the conditions of Assertion 4.3 are satisfied for the function $g(n)=n^l, l\geq 0$.

$g(n)$ - increases, $\lim_{n\to\infty} g(n)=\infty$, $\text{g}'(n)$ - continuous function,

$$\lim_{n\to\infty}\frac{g(n)}{\ln(n)g'(n)}=\lim_{n\to\infty}\frac{n^l}{\ln(n)\,\mathrm{l(n)}^{l-1}}=\infty\neq 0.$$

Let us show that the conditions of Assertion 4.3 are not satisfied for the function $g(n)=2^n$.

$g(n)$ - increases, $\lim_{n\to\infty} g(n)=\infty$, $\text{g}'(n)$ - continuous function,

$$\lim_{n\to\infty}\frac{g(n)}{\ln(n)g'(n)}=\lim_{n\to\infty}\frac{2^n}{\log n 2^n \log 2}=0.$$

Now let's consider a necessary condition for the above asymptotics to hold.

We will use the notation of Assertion 4.1.

Assertion 4.4

Let:

$$\lim_{n\to\infty}\frac{\sum_{m=1}^{n}a_m g(m)}{\sum_{m=1}^{n}b_m g(m)}=1$$

Then it is fulfilled, when $p\to\infty$ ( $p$ - a prime number ):

$$\left|\frac{g(p)}{\sum_{m=1}^{p}b_m g(m)}\right|\to 0\,. \tag{4.10}$$

Proof

If:

$$\lim_{n\to\infty}\frac{\sum_{m=1}^{n}a_m g(m)}{\sum_{m=1}^{n}b_m g(m)}=1$$

Then it executes:

$$\lim_{n\to\infty}\left|\frac{\sum_{m=1}^{n}a_m g(m)}{\sum_{m=1}^{n}b_m g(m)}-\frac{\sum_{m=1}^{n-1}a_m g(m)}{\sum_{m=1}^{n-1}b_m g(m)}\right|=0\,. \tag{4.11}$$

Let's take $n=p$, where $p$ is a prime number. Then having in mind (3.11):

$$\lim_{n\to\infty}\left|\frac{\sum_{m=1}^{p}a_m g(m)}{\sum_{m=1}^{p}b_m g(m)}-\frac{\sum_{m=1}^{p-1}a_m g(m)}{\sum_{m=1}^{p-1}b_m g(m)}\right|=0 \tag{4.12}$$

Let us denote the value $a_m g(m)=a_p g(p)$ when $m=p$.

Let's use that $a_p=1$ and we get:

$$\sum_{m=1}^{p}a_m\,\mathrm{g}(m)=a_p g(p)+\sum_{m=1}^{p-1}a_m g(m)=g(p)+\sum_{m=1}^{p-1}a_m g(m)\,. \tag{4.13}$$

We transform (4.12) and, taking into account (4.13), we obtain:

$$\left|\frac{\sum_{k=1}^{p} a_m g(m)}{\sum_{k=1}^{p} b_m g(m)} - \frac{\sum_{k=1}^{p-1} a_m g(m)}{\sum_{k=1}^{p-1} b_m g(m)}\right| = \left|\frac{a_p g(p)\sum_{k=1}^{p-1} b_m g(m) - b_p g(p)\sum_{k=1}^{p-1} a_m g(m)}{\sum_{m=1}^{p} b_m g(m)\sum_{k=1}^{p-1} b_m g(m)}\right| = |f(p)|\left|\frac{1 - b_p \frac{\sum_{m=1}^{p-1} a_m g(m)}{\sum_{m=1}^{p-1} b_m g(m)}}{\sum_{m=1}^{p} b_m g(m)}\right|$$

Having in mind that at $p \to \infty$ the value $b_p \to 0$ and $\lim_{n\to\infty}\frac{\sum_{m=1}^{n} a_m g(m)}{\sum_{m=1}^{n} b_m g(m)} = 1$, then it is performed:

$$\left|\frac{g(p)}{\sum_{m=1}^{p} b_m g(m)}\right| \to 0,$$

which corresponds to (4.10).

Let's look at examples of fulfilling the necessary condition:

1. $g(p) = \log p$. Since $\sum_{m=1}^{p} b_m g(p) = d^*(f)\sum_{m=2}^{p}\frac{\ln m}{\ln m} = d^*(f)p$, therefore

$$\left|\frac{g(p)}{\sum_{m=1}^{p} b_m f(p)}\right| = \left|\frac{\ln p}{d^*(f)p}\right| \to 0 \text{ at } p \to \infty.$$

2. $g(p) = p^l, l > 0$. Since $\sum_{m=1}^{p} b_m g(p) = d^*(f)\sum_{m=2}^{p}\frac{m^l}{\ln m} = \frac{d^*(f)p^{l+1}}{\ln(\mathrm{p})}(1+o(1))$.

Therefore $\left|\frac{g(p)}{\sum_{m=1}^{p} b_m g(p)}\right| = \left|\frac{p^l}{\frac{d^*(f)p^{l+1}}{\ln(\mathrm{p})}(1+o(1))}\right| \to 0$, at $p \to \infty$.

Now let's look at an example of not meeting the necessary condition:

3. $g(p)=\dfrac{1}{p^2}$. Since $\sum_{m=1}^{p} b_m g(p) = d^*(f)\sum_{m=2}^{p}\dfrac{1}{m^2\ln m} = C$, where $C$ is a constant.

That's why $\lim_{p\to\infty} \left| \dfrac{g(p)}{\sum_{m=1}^{p} b_m g(p)} \right| = \lim_{p\to\infty} \left| \dfrac{\frac{1}{p^2}}{C} \right| \neq 0$.

## 4. CONCLUSION AND SUGGESTIONS FOR FURTHER WORK

The next article will continue to study the behavior of some sums of functions of primes.

## 5. ACKNOWLEDGEMENTS

Thanks to everyone who has contributed to the discussion of this paper. I am grateful to everyone who expressed their suggestions and comments in the course of this work.